\documentclass[12pt]{article}
\usepackage{amssymb,amsbsy,amsmath,amsfonts,amssymb}
\usepackage{latexsym,euscript,exscale}

\usepackage{times}

\newcommand{\pf}{

\smallskip

\noindent {\it Proof : }}

\newcommand {\N}{\mathbb N}

\newcommand {\R}{\mathbb R}

\newcommand {\C}{\mathbb C}
\newcommand {\HH}{\mathbb H}

\newcommand {\DD}{{\mathbb D}_2}

\newcommand{\pff}{$\hfill \square$
\smallskip}

\newcommand{\norm}[1]{\ensuremath{\left\|#1\right\|}}

\newtheorem{prop}{Proposition}
\newtheorem{lemm}[prop]{Lemma}
\newtheorem{theo}[prop]{Theorem}
\newtheorem{ques}[prop]{Question}
\newtheorem{coro}[prop]{Corollary}

\newtheorem{defi}[prop]{Definition}

\title{Even infinite dimensional real Banach spaces}

\author{Valentin Ferenczi\footnote{This article was written during a visit of the first author at the University of S\~ao Paulo financed by the Centre National de la Recherche Scientifique.} \ and El\'oi  Medina Galego}

\date{ }

\begin{document}

\maketitle

\begin{abstract} This article is a continuation of a paper of the first author \cite{F} about complex structures on real Banach spaces. We define a notion of even infinite dimensional real Banach space, and prove that there exist even spaces,  including HI or unconditional examples from \cite{F} and $C(K)$ examples due to Plebanek \cite{P}. We extend results of \cite{F} relating the set of complex structures up to isomorphism on a real space to a group associated to inessential operators on that space, and give characterizations of even spaces in terms of this group.  We also generalize results of \cite{F} about totally incomparable complex structures to essentially incomparable complex structures, while showing that the complex version of a space defined by S. Argyros and A. Manoussakis \cite{AM} provide examples of essentially incomparable complex structures which are not totally incomparable.
\end{abstract}

\section{Introduction}

Any complex Banach space is also a real Banach space. Conversely, the
linear structure on a real Banach space $X$
may be induced by a $\C$-linear structure; the corresponding complex Banach space is
said to be a {\em complex structure} on $X$.

The existing theory of complex structure is up to isomorphism. In this setting, a complex structure  on a real Banach space $X$ is a complex space which is $\R$-linearly isomorphic to $X$. Any complex structure up to isomorphism is associated to an $\R$-linear isomorphism $I$ on $X$ such that $I^2=-Id$. Conversely, for any such isomorphism $I$, an associated complex structure may be defined by the law 
 $$\forall \lambda,\mu \in \R, (\lambda+i\mu).x=(\lambda Id+\mu
I)(x),$$
and  the equivalent norm
$$|||x|||=\sup_{\theta \in \R}\norm{\cos\theta x+\sin\theta Ix}.$$

Isomorphic theory of complex structure addresses questions of existence, uni\-queness, and when there is more than one complex structure, the possible structure of the set of complex structures up to isomorphism (for example in terms of cardinality).

It is well-known that complex structures do not always exist (up to isomorphism) on a Banach space.
The HI space of Gowers and Maurey \cite{GM1} is a good example of this, or more generally any space with the $\lambda Id+S$ property (i.e. every operators is a strictly singular perturbation of a multiple of the identity); and note that since this property passes to hyperplanes, complex structures neither exist on finite codimensional subspaces of such a space (relate this observation to the forthcoming Question \ref{q}).

By the examples of \cite{B} and \cite{K} of complex spaces
not isomorphic to their conjugates, there exists real spaces with at least two complex structures up to isomorphism, and the examples of \cite{B} and \cite{An} (which are separable) actually admit 
 a continuum of complex structures.
In \cite{F} the first author showed that for each $n \geq 2$ there exists a space with exactly $n$ complex structures. He also gave various examples of spaces with unique complex structure up to isomorphism and different from the classical example of $\ell_2$, including a HI example and a space with an unconditional basis.

\

A fundamental tool
in \cite{F} is an identification of isomorphism classes of complex structures on a space $X$ with conjugation classes in some group associated to strictly singular operators on $X$. It remained open whether the associated map was bijective. In this paper, we show that it is not bijective in general, but that there actually exists a natural bijection between isomorphism classes of complex structure on a space $X$ and on its hyperplanes with conjugation classes in the group associated to strictly singular operators in $X$, Theorem \ref{resume}. We also note that this holds as well when one replaces strictly singular operators by any Fredholm ideal in ${\cal L}(X)$,  such as $In(X)$ the ideal of inessential operators.

More precisely, it turns out to be fundamental to determine when a given operator of square $-Id$ modulo inessential operators lifts to an operator of square $-Id$ on the space. While the answer is always positive in the complex case, Proposition \ref{com}, it turns out that on a real space, such an operator lifts either to an operator on $X$ of square $-Id$
or, in a sense made precise in Lemma \ref{hyperp}, to an operator
on an hyperplane of $X$ of square $-Id$, Proposition \ref{disjonction}; furthermore the two cases are exclusive, Proposition \ref{exclusion}. This implies that operators of square $-Id$ modulo inessential operators characterize complex structures on $X$ and on its hyperplanes.

This leads us to define a notion of even and odd real Banach space extending the classical notion for finite dimensional spaces. Even spaces are those spaces which admit complex structure but whose hyperplanes do not. Odd spaces are the hyperplanes of even spaces.

We provide characterizations of even and odd spaces in terms of the previously mentioned groups associated to inessential operators on $X$ and on its hyperplanes, Corollary \ref{oddeven} or more precisely, in terms of lifting properties of operators of square $-Id$ modulo inessential operators, Proposition \ref{equivalence}. 
We prove that there exist even infinite dimensional Banach spaces,
using various examples from \cite{F}, including a HI and an unconditional example, Theorem \ref{existseven}. Moreover we use  spaces  constructed in \cite{P} to give examples
of even and odd spaces of the form $C(K)$, Theorem \ref{CK}. We also show that the direct sum of essentially incomparable infinite dimensional spaces is even whenever both spaces are even, Proposition \ref{evenpluseven}.

Finally we extend and simplify the proof of some results of \cite{F} about totally incomparable complex structures by showing that essentially incomparable complex structures are necessarily conjugate, Theorem \ref{ttt} and Corollary \ref{ccc}.
We also show that the complex version of a space built by S. Argyros and A. Manoussakis \cite{AM} provides examples of complex structures on a space which  
are essentially incomparable yet not totally incomparable, Proposition \ref{exam}.

\

\section{Parity of infinite dimensional spaces}

It may be natural to think that a real infinite dimensional
space of the form $X \oplus X$ should be considered to be even.
 This seems to be restrictive however, as we should consider as candidates
 for spaces with even dimension the spaces which admit a complex structure,
 and there are  spaces with complex structure which are not isomorphic to a
cartesian square (actually, not even decomposable,
 by the examples of \cite{F}).
Another problem is that we would wish the hyperplanes of a space with even infinite dimension not to share the same property. In other words, parity should imply a structural difference between the whole space and its hyperplanes.
This suggests the following definition, which obviously generalizes the case of finite dimensional spaces, and will be our guideline for this section.

\begin{defi}\label{even}
A real Banach space is {\em even} if it admits a complex structure but its hyperplanes do not admit a complex structure. It is {\em odd} if its hyperplanes are even.
\end{defi}

Equivalently a Banach space is odd  if it does not have a complex structure
but its hyperplanes do, and clearly $2$-codimensional subspaces of an even
(resp. odd) space are  even (resp. odd).

\

The following crucial fact will be used repeatedly without explicit reference: two complex structures $X^I$ and $X^J$ are isomorphic if and only if $I$ and $J$ are conjugate, i.e. there exists an isomorphism $P$ on $X$ such that $J=PIP^{-1}$ ($P$ is then $\C$-linear from $X^I$ onto $X^J$). There is therefore a natural correspondance between isomorphism classes of complex structure on $X$ and conjugacy classes of elements of square $-Id$ in ${\cal GL}(X)$, and we shall sometimes identify the two sets. 

\

 Our first results are 
improvements of  results from \cite{F}. We recall that an operator $T: Y \rightarrow Z$ is Fredholm if its kernel is finite dimensional and its image is finite codimensional,
in which case the Fredholm index of $T$ is defined by
 $$i(T)=\dim(Ker(T))-\dim(Z/TY).$$ 
We shall use the easy facts that an operator $T$ is Fredholm if and only if $T^2$ is Fredholm (with $i(T^2)=2i(T)$), and that
a $\C$-linear operator is Fredholm if and only if it is Fredholm as an $\R$-linear operator (and the corresponding indices are related by $i_{\R}(T)=2i_{\C}(T)$).

\

 A closed two-sided
 ideal ${\cal U}(X)$ in ${\cal L}(X)$ is a {\em Fredholm ideal} when an
 operator $T \in {\cal L}(X)$ is Fredholm if and only if
 the corresponding class is invertible in ${\cal L}(X)/{\cal U}(X)$.
 It follows from well-known results in Fredholm theory that ${\cal U}(X)$
 is contained in the ideal  ${\cal I} n(X)$ of inessential operators, i.e. operators $S$ such that
$Id_X-VS$ is a Fredholm operator for all operators $V \in {\cal L}(X)$,
 see for example \cite{Gon}. Note that by continuity of the Fredholm index, $Id_X-VS$ is necessarily Fredholm with index $0$ when $S$ is inessential, and $V \in {\cal L}(X)$.

\begin{lemm}\label{perturbation} Let $X$ be an infinite dimensional real Banach space.
 Let $I \in {\cal L}(X)$ satisfy $I^2=-Id$,
 and let $S$ be inessential
 such that $(I+S)^2=-Id$. Then $I$ and $I+S$ are conjugate, or equivalently, the complex structures $X^I$ and $X^{I+S}$ associated to $I$ and $I+S$
respectively are isomorphic.
\end{lemm}

\pf The map $2I+S$ is immediately seen to be $\C$-linear from $X^I$ into
$X^{I+S}$. Furthermore, it is an inessential perturbation of an isomorphism,
 and therefore Fredholm with index $0$ as an $\R$-linear operator on $X$.
So it is  also Fredholm with index $0$ as a $\C$-linear map, which implies that $X^I$ and $X^{I+S}$ are $\C$-linearly isomorphic. \pff

\

Let $X$ be a Banach space. The set ${\cal I}(X)$ denotes the set of operators
on $X$ of square $-Id$.
 Let ${\cal U}(X)$ be a Fredholm ideal
 in ${\cal L}(X)$ and $\pi_{\cal U}$ (or $\pi$) denote the quotient map
from
${\cal L}(X)$ onto ${\cal L}(X)/{\cal U}(X)$.  Let
$({\cal L}(X)/{\cal U}(X))_0$ denote the group $\pi_{\cal U}({\cal GL}(X))$, and
 $\tilde{{\cal I}}(X)$ denote the set of elements of $({\cal L}(X)/{\cal U}(X))_0$ whose square is
equal to $-\pi_{\cal U}(Id)$.

\begin{lemm}\label{principe} Let $X$ be an infinite dimensional real Banach space and ${\cal U}(X)$ be a Fredholm ideal in ${\cal L}(X)$. Then the quotient map $\pi_{\cal U}$ induces an injective map
$\tilde{\pi_{\cal U}}$ from the set of ${\cal GL}(X)$-conjugation classes on
${\cal I}(X)$ (and therefore from the set of isomorphism classes of complex structures on $X$) into the set of $({\cal L}(X)/{\cal U}(X))_0$-conjugation classes on
$\tilde{{\cal I}}(X)$. 
\end{lemm}

\pf For any operator $T$ on $X$, we write $\tilde{T}=\pi(T)$. Let $I$ and $T$ be operators in ${\cal I}(X)$. If
$\alpha$ is a $\C$-linear isomorphism from
$X^I$ onto $X^T$, then the $\C$-linearity means that $\alpha I=T \alpha$.
Therefore $\tilde{\alpha}\tilde{I}=\tilde{T}\tilde{\alpha}$, and $\tilde{I}$ and $\tilde{T}$ satisfy a conjugation
relation. Conversely, if $\tilde{I}=\tilde{\alpha}^{-1} \tilde{T} \tilde{\alpha}$ for some $\alpha \in {\cal GL}(X)$,
then  $\alpha^{-1}T\alpha=I+S$, where $S$ belongs to ${\cal U}(X)$ and is therefore inessential .
Note that
 $(I+S)^2=-Id$, and since $T\alpha=\alpha(I+S)$, $\alpha$ is a $\C$-linear isomorphism from 
$X^{I+S}$ onto $X^T$.
By Lemma \ref{perturbation}, it follows that $X^I$ and $X^T$ are isomorphic.
This proves that $\tilde{\pi}$ is well-defined and injective.
 \pff

\

We shall now discuss when the above induced map is actually a bijection.
 This is equivalent to saying that ${\cal U}(X)$ has the following lifting property.

\begin{defi}
Let $X$ be an infinite dimensional Banach space.
The Fredholm ideal ${\cal U}(X)$ is said to have
 the {\em lifting property} if any $\alpha$
 in $({\cal L}(X)/{\cal U}(X))_0$ satisfying $\alpha^2=-1$
 is the image under $\pi_{\cal U}$ of an operator $T$ such that $T^2=-Id$.
\end{defi}

The following was essentially observed in \cite{F}. 

\begin{lemm} \label{subalgebra}
Let $X$ be an infinite dimensional Banach space and let ${\cal U}(X)$ be a Fredholm ideal in ${\cal L}(X)$. If  ${\cal U}(X)$ admits a supplement in ${\cal L}(X)$ which is a sub\-algebra of ${\cal L}(X)$,
then ${\cal U}(X)$ has the lifting property.
\end{lemm}

\pf If ${\cal H}(X)$ is a subalgebra of ${\cal L}(X)$ which supplements ${\cal S}(X)$, then let $T \in {\cal L}(X)$ be
such that $\tilde{T}^2=-\tilde{Id}$; we may assume that $T$ (and therefore $T^2$) belongs to ${\cal H}(X)$. Then
since
$T^2+Id$ is in ${\cal U}(X) \cap {\cal H}(X)$,
 $T^2$ must be equal to $-Id$. Any class $\tilde{T} \in \tilde{\cal I}(X)$ may  therefore be
lifted to an element of ${\cal I}(X)$.\pff

\

We shall now prove that although any Fredholm ideal in a complex space has the
lifting property, this is not necessarily true in the real case. 
The proof of the complex case is essentially the same as the similar  classical result concerning projections (see e.g. \cite{Gon}), and could be deduced directly from it using the fact that an operator $A$ satisfies $A^2=-Id$ if and only if $\frac{1}{2}(Id-iA)$ is a projection. We shall however give a direct proof of this result for the sake of completeness.

\begin{prop}\label{com} Let $X$ be an infinite dimensional complex Banach space and let ${\cal U}(X)$ be a Fredholm ideal in ${\cal L}(X)$. Then every  element $a \in {\cal L}(X)/{\cal U}(X)$ with $a^2=-1$ is  image under the quotient map of some $A \in {\cal L}(X)$ with $A^2=-Id$.
\end{prop} 
 
 \pf Recall that  $\pi: {\cal L}(X) \to {\cal L}(X)/{\cal U}(X)$ denote the quotient map. We choose $B \in {\cal L}(X)$ such that $\pi(B)=a$. So $\pi(B^2)=-1$ and therefore there exists $S \in {\cal U}(X) \subset {\cal I}n(X)$ such that $B^{2}=-Id+S$.
 
 Since the spectrum $\sigma(-Id+S)$ of $-Id+S$ is countable and its possible
 limit point  is -1, it follows by the spectral mapping  theorem (\cite{DS}, Theorem VII.3.11) that the spectrum $\sigma(B)$ is also countable and its possible limit points are $-i$ and $i$.
 
 Take a simple closed curve $\Gamma$ in $\C \setminus \sigma(B)$ such that $i$ is enclosed by $\Gamma$ and $-i$ is not enclosed by $\Gamma$.  Define the operator
 $$P= \int_{\Gamma} (\lambda I - B)^{-1} d{\lambda}.$$
 By \cite{DS}, Theorem VII.3.10, $P$ is a projection. Moreover, according to the continuity of $\pi$,
  $$\pi(P)= \int_{\Gamma} (\lambda I - a)^{-1} d{\lambda}. \eqno(1)$$
  On the other hand, it is easy to check that
   $$(\lambda I -a)^{-1}= \frac{\lambda I + a}{{\lambda}^{2}+1}$$
   Thus (1) implies that 
   $$\pi(P)=\frac{i+a}{2i},$$
 and hence putting $A=2iP-iI$, we see that $\pi(A)=a$ and $A{^2}=-I$. \pff

\

Such a proof is not possible in the real case essentially because there is no
formula with real coefficients linking projections and operators of square $-Id$. Actually it is known that if $X$ is real  every element $p$ belonging to the quotient by a Fredholm ideal
and satisfying $p^2=p$ may be lifted to a projection (see \cite{Gon}), and the
proof uses the complexification of $X$ and a curve with well-chosen symmetry so that the complex
projection produced by the proof in the complex case is induced by a real
projection which will answer the question by the positive. However there is no choice of curve
such that the complex map $A$ of square $-Id$  obtained in the above proof applied in the complexification of $X$ is induced by a real
operator. Actually the result is simply false in the real case:

\begin{theo} Let $X$ be a real infinite dimensional Banach space whose hyperplanes admit a complex
  structure. Then no Fredholm ideal in ${\cal L}(X)$  has the lifting property. \end{theo}

This theorem is a consequence of Proposition \ref{exclusion}.
We need to recall that the complexification $\hat X$ of a real Banach space $X$ (see, for example, \cite{LT}, page 81) is defined as the space 
 $\hat X= \{x+iy: x, y \in X \},$ which is the space $X \oplus X$ with the canonical complex structure associated to $J$ defined on $X \oplus X$ by $J(x,y)=(-y,x)$.
 Let $A, B \in L(X)$. Then 
 $$(A+iB)(x+iy):=Ax-By+(Ay+Bx)$$
 defines an operator $A+iB \in L(\hat X)$ that satisfies 
 $\max \{ \|A \|, \|B \| \} \leq \|A+iB \| \leq 2^{1/2} (\|A \|+ \|B\|).$
 Conversely, given $T \in L(\hat X)$, if we put $T(x+i0):=Ax+iBx$, then we obtain $A, B \in L(X)$ such that $T=A+iB$.
 We write $\hat{T}=T+i0$ for $T \in {\cal L}(X)$, and say that such the operator $\hat{T}$ is induced by the real operator $T$.

\begin{prop}\label{exclusion} Let $Y$ be an infinite dimensional real Banach space and $J \in {\cal L}(Y)$ such that $J^2=-Id$. Let $A$ be defined on
$X=\R \oplus Y$  by the matrix 
$\begin{pmatrix} 1 & 0 \\ 0 & J \end{pmatrix}$. Then $A^2$ is the sum of $-Id$ and of a rank $1$ operator, but there is no inessential operator $S$ on $X$ such that $(A+S)^2=-Id$.
\end{prop}

\pf Assume $(A+S)^2=-Id$ for some inessential operator $S$. Passing to the complexification $\hat{X}$ of $X$, we obtain
$(\hat{A}+\hat{S})^2=-\hat{Id}$. The map from $[0,1]$ into ${\cal L}(\hat{X})$ defined by $T_{\mu}=\hat{A}+\mu \hat{S}$ is polynomial, moreover by spectral properties of inessential operators and the spectral theorem, the spectrum 
$Sp(T_\mu)$ is, with the possible exception of $i$ and $-i$, a countable set of isolated points, which are eigenvalues with associated spectral projections, denoted 
$E(\lambda,T_\mu)$ for each $\lambda \in Sp(T_\mu)$, of finite rank. Furthermore the complex operator $T_\mu$ is induced by the real operator $A+\mu S$, therefore $Sp(T_\mu)$ is symmetric with respect to the real line.

Let $n(\mu)=\sum_{\lambda \in \R \cap Sp(T_{\mu})}rk(E(\lambda,T_\mu))$ and
let $I_1=\{\mu \in [0,1]: n(\mu) \ is \ odd\}$,
$I_0=\{\mu \in [0,1]: n(\mu) \ is \ even\}$.

Observe that $0 \in I_1$, since $\hat{A}$ is defined on $\hat{X}=\hat{\R} \oplus \hat{Y}$ by the matrix
$\begin{pmatrix} 1 & 0 \\ 0 & \hat{J} \end{pmatrix}$, and therefore has unique real eigenvalue $1$, with associated spectral projection of dimension $1$. On the other hand, since $T(1)^2=(\hat{A}+\hat{S})^2=-\hat{Id}$, $T(1)$ does not admit real eigenvalues and therefore $1 \in I_0$.

We now pick some $\mu_1 \in I_1$. Let $U$ be an open set containing the real line, symmetric with respect to it, and such that $\overline{U} \cap Sp(T_{\mu_1}) \subset \R$.
Then the spectral projection $E(\mu):=E(T_\mu,U \cap Sp(T_\mu))$ is an analytic (\cite{DS} Lemma VII.6.6) projection valued function defined for all $\mu$ such that $|\mu-\mu_1|<\gamma$ for some small enough $\gamma>0$, and for which
$E(\mu_1)$ has rank $n(\mu_1)$. Therefore by \cite{DS} Lemma VII.6.8, the dimension $\sum_{\lambda \in U \cap Sp(T_\mu)} rk(E(\lambda,T_\mu))$ of the image of $E(\mu)$ is also $n(\mu_1)$, for $\mu$ in a small enough open set $V$ around $\mu_1$. 
By symmetry of $U$ and of $Sp(T_\mu)$ with respect to the real line (with preservation of the ranks of the associated spectral projections),
$$n(\mu_1)=\sum_{\lambda \in \R \cap Sp(T_\mu)} rk (E(\lambda, T_\mu))+2 \sum_{\lambda \in Sp(T_\mu), Im(\lambda)>0} rk (E(\lambda, T_\mu)),$$
when $\mu \in V$.
So $n(\mu)=\sum_{\lambda \in \R \cap Sp(T_\mu)} rk (E(\lambda,T_\mu))$ is odd whenever $\mu$ is in the neighborhood $V$ of $\mu_1$.

We have therefore proved that $I_1$ is open. In the same way, $I_0$ is also open (in the special case when $n(\mu_0)=0$, then $E(\mu_0)=0$ and so, $E(\mu)=0$ and therefore
$n(\mu)=0$ in a neighborhood of $\mu_0$). In conclusion, the sets $I_0$ and $I_1$ are open, non-empty, and partition $[0,1]$, a contradiction. \pff

\

The obstruction for the lifting property is therefore that a complex structure
on a hyperplane of a space $X$ does not correspond to a complex structure on $X$, although it does induce elements of square $-1$ in 
$({\cal L}(X)/{\cal S}(X))_0$, as explicited by the following result.

\begin{lemm}\label{hyperp} Let $Y$ be an infinite dimensional Banach space and $X=Y \oplus \R$. Let ${\cal U}(X)$ be a Fredholm ideal. Then the map $\pi_{\cal U}^{\prime}$
(or $\pi'$) from ${\cal GL}(Y)$ into $({\cal L}(X)/{\cal U}(X))_0$ defined by
$\pi_{\cal U}^{\prime}(A)=\pi_{\cal U}(\begin{pmatrix} 1 & 0 \\ 0 & A \end{pmatrix})$
maps ${\cal I}(Y)$ into $\tilde{{\cal I}}(X)$ and induces an injection from the set of conjugation classes on ${\cal I}(Y)$ (and therefore from the set of isomorphism classes of complex structures on $Y$) into the set of conjugation classes on $\tilde{{\cal I}}(X)$. \end{lemm}

\pf First note that if the conjugation relation $J=PKP^{-1}$ is satisfied in ${\cal GL}(Y)$ then 
$$ \begin{pmatrix} 1 & 0 \\ 0 & J \end{pmatrix}=
\begin{pmatrix} 1 & 0 \\ 0 & P \end{pmatrix}
\begin{pmatrix} 1 & 0 \\ 0 & K \end{pmatrix}
\begin{pmatrix} 1 & 0 \\ 0 & P^{-1} \end{pmatrix}$$
is satisfied in ${\cal GL}(X)$, which provides a conjugation relation
in $({\cal L}(X)/{\cal U}(X))_0$. Conversely if 
$$ \widetilde{\begin{pmatrix} 1 & 0 \\ 0 & J \end{pmatrix}}= 
\tilde{P}
\widetilde{\begin{pmatrix} 1 & 0 \\ 0 & K \end{pmatrix}}
\tilde{P}^{-1},$$
write $$P=\begin{pmatrix} a & b^* \\ c & D \end{pmatrix}.$$
Then since $P$ is Fredholm with index $0$, $P_{|Y}$ is Fredholm with index $-1$ as an operator of ${\cal L}(Y,X)$, and  $D=P_{|Y}-b^*$ is also Fredholm with index $-1$ in
${\cal L}(Y,X)$, and with index $0$ in ${\cal L}(Y)$. Therefore some finite rank perturbation $D'$ of $D$ is an isomorphism. Furthermore it is easy to deduce from the conjugation relation that 
$J=D'K{D'}^{-1}+S$, where $S \in {\cal U}(Y)$. Therefore
$J-S=D'K{D'}^{-1}$, i.e. $J-S$ and $K$ are conjugate, and by
Lemma \ref{principe}, it follows that $Y^J$ is isomorphic to $Y^K$. \pff

\

The previous lemma does not mean that the sets of isomorphism classes of complex structures on $X$ and on hyperplanes of $X$ are disjoint (only the correspoding images by $\tilde{\pi}$, and $\tilde{\pi'}$ respectively, are). Actually the sets of isomorphism classes on $X$ and on its hyperplanes can either be equal, when $X$ is isomorphic to its hyperplanes, or disjoint, when it is not. And we have:

\begin{prop}\label{disjonction} Let $Y$ be an infinite dimensional Banach space and let $X=\R \oplus Y$.
Let ${\cal U}(X)$ be a Fredholm ideal in ${\cal L}(X)$.
 Let $A \in GL(X)$ and assume $A^2=-Id+S$, $S \in {\cal U}(X)$ (i.e. $\tilde{A} \in \tilde{{\cal I}}(X)$). Then there exists 
$s \in {\cal U}(X)$ such that either $(A+s)^2=-Id$ or
$(A+s)^2=\begin{pmatrix} 1 & 0 \\ 0 & J \end{pmatrix}$ where
$J \in {\cal L}(Y)$ satisfies $J^2=-Id$. \end{prop}

Before the proof, let us observe that by Proposition \ref{exclusion}, only one of the two alternatives of the conclusion can hold for a given $A$ such that
$\tilde{A} \in \tilde{\cal I}(X)$.

\

\pf Passing to the complexification $\hat{X}$ of $X$, we have that $\hat{A}^2=-\hat{Id}+\hat{S}$, and $\hat{S}$ is inessential. Now let $\Gamma$ be a rectangular curve with horizontal and vertical edges, symmetric with respect to the horizontal axis, included in the open unit disk, and such that
$\Gamma \cap Sp(\hat{S}) = \emptyset$ and let $U$ be the  interior of the domain delimited by $\Gamma$, $V$ be the interior of the complement of this domain.
Let $\hat{P}$ be the spectral projection associated to
$Sp(\hat{S}) \cap U$; since $\Gamma$ is rectangular and symmetric with respect to the real axis, it is classical and easy to see that $\hat{P}$ is indeed induced by a real operator $P$ on $X$, see e.g. \cite{Gon2} where this principle is used.  Let also
$\hat{Q}$ be the spectral projection associated to
$Sp(\hat{S}) \cap V$. 

Then $\hat{S}=\hat{S}\hat{P}+\hat{S}\hat{Q}$. The operator
$\hat{S}\hat{P}$ has spectral radius strictly smaller than $1$, therefore the series $\sum_{n \geq 1}
b_n (\hat{S}\hat{P})^n$ converges to an operator $\hat{s}$, where
$\sum_{n \geq 1} b_n z^n=-1+(1-z)^{-1/2}$ for all $|z|<1$, and since the $b_n$'s are reals,  it is indeed induced by a real operator $s=\sum_{n \geq 1} b_n (SP)^n$ in ${\cal U}(X)$.
We observe that
$$(\hat{P}+\hat{s})^2=\hat{P}(\hat{Id}+\hat{s})^2=\hat{P}
(\hat{Id}-\hat{S})^{-1},$$
therefore
$$(\hat{A}\hat{P}+\hat{A}\hat{s})^2=-\hat{P}.$$

Assume now that $Q$ has even rank, then there exists a finite rank operator $F$ on $QX$ such that $F^2=-Id_{|QX}$. Let then
$s'=FQ-AQ$, therefore $\hat{A}\hat{Q}+\hat{s'}=\hat{F}\hat{Q}$ and $(\hat{A}\hat{Q}+\hat{s'})^2=-\hat{Q}$.
If we then let $v=s+As'$, we deduce that
$$(\hat{A}+\hat{v})^2=(\hat{A}\hat{P}+\hat{A}\hat{s}+\hat{A}\hat{Q}+\hat{s'})^2=-\hat{P}-\hat{Q}=-\hat{Id},$$
and therefore $(A+v)^2=-Id$, with $v \in {\cal U}(X)$.

If $Q$ has odd rank then there exists a finite rank operator $F$ on $QX$
 such that $F=\begin{pmatrix} 1 & 0 \\ 0 & j \end{pmatrix}$, with $j^2=-Id$,
  in an appropriate decomposition of $QX$. Defining $v \in {\cal U}(X)$ in the same way as above, we obtain
 that $A+v$ may be written $$A+v= \begin{pmatrix} 1 & 0 \\ 0 & J \end{pmatrix},$$ 
corresponding to some decomposition $R'\oplus Y'$ of $X$ where $R'$
 is $1$-dimensional.
If $Y'=Y$ then we may clearly find some rank $1$ perturbation $f$ such that $A+v+f$ may be written
$$A+v+f= \begin{pmatrix} 1 & 0 \\ 0 & J \end{pmatrix},$$
corresponding to the original decomposition $\R \oplus Y$ of $X$.
If $Y'\neq Y$ then  we consider the space $Z=Y'\cap Y \cap JY$, which is stable by
$J$. If $Z$ has codimension $3$ then we may decompose $Y=G \oplus Z$,
where $G$ has dimension $2$, and by using a operator $k$ of square $-Id_{G}$ on
$G$, we may find a rank $3$ perturbation $f$ such that
$$A+v+f= \begin{pmatrix} 1 & 0 \\ 0 & J' \end{pmatrix},$$
in the original decomposition $\R \oplus Y$ of $X$, and with
$J'_{|Z}=J_{|Z}$ and $J'_{|G}=k$, so that $J^{\prime 2}=-Id_{|Y}$.
If finally $Z$ has codimension $2$ then easily some rank $2$
perturbation of $A+v$ on $X$ has square $-Id$, but because of the decomposition
of $A+v$ on $R'\oplus Y'$, this would contradict Proposition \ref{exclusion}, so this case is not possible. 
\pff

\

We sum up the results of Proposition \ref{exclusion}, Lemma
\ref{hyperp} and Proposition \ref{disjonction} in the next theorem:

\begin{theo}\label{resume} Let $Y$ be an infinite dimensional real Banach space and let $X=Y \oplus \R$. Let ${\cal U}(X)$ be a Fredholm ideal in ${\cal L}(X)$. Then there exists a partition
$\{\tilde{{\cal I}_0}(X),\tilde{{\cal I}_1}(X)\}$ of
$\tilde{{\cal I}}(X)$ such that
$\pi_{\cal U}$ induces a bijection from the set of conjugation classes on ${\cal I}(X)$ (and therefore from the set of isomorphism classes of complex structures on $X$) onto the set of conjugation classes on
$\tilde{{\cal I}_0}(X)$ and such that $\pi_{\cal U}^{\prime}$ induces a bijection from the set of conjugation classes on ${\cal I}(Y)$ (and therefore from the set of isomorphism classes of complex structures on $Y$) onto the set of conjugation classes on
 $\tilde{{\cal I}_1}(X)$. \end{theo}

\begin{coro}\label{oddeven}
An infinite dimensional real Banach space is even if and only if
$\tilde{\cal I}(X) =\tilde{{\cal I}_0}(X) \neq \emptyset$, and odd
if and only if $\tilde{\cal I}(X)=\tilde{{\cal I}_1}(X) \neq \emptyset$.
When $\tilde{\cal I}(X)$ is a singleton then $X$ is either even or odd.
\end{coro}

The next proposition sums up when a real Banach space has the lifting property.

\begin{prop}\label{equivalence} Let $X$ be an infinite dimensional real Banach space. Then the following are equivalent:
 
\begin{itemize}
\item i) any Fredholm ideal in ${\cal L}(X)$ has the lifting property.

\item ii) some Fredholm ideal in ${\cal L}(X)$ has the lifting property.

\item iii) the hyperplanes of $X$ do not admit complex structure.

\item iv) for any Fredholm ideal ${\cal U}(X)$,
the map $\pi_{\cal U}$ induces a bijection from the set of isomorphism classes of
complex structures on $X$ onto the set of conjugation classes on
$\tilde{{\cal I}}_{\cal U}(X)$.

\item v) for some Fredholm ideal ${\cal U}(X)$,
the map $\pi_{\cal U}$ induces a bijection from the set of isomorphism classes of
complex structures on $X$ onto the set of conjugation classes on
$\tilde{{\cal I}}_{\cal U}(X)$.
\end{itemize}
\end{prop}

\pf It is clear by definition that $iv) \Leftrightarrow  i)$ and 
$v) \Leftrightarrow ii)$. Then $i) \Rightarrow ii)$ is obvious,
$ii) \Rightarrow iii)$ due to Proposition \ref{exclusion}, and $iii) \Rightarrow i)$ by Proposition \ref{disjonction}.
\pff

\

We use this proposition to solve an open question from the first author, which in our formulation asked whether there existed infinite dimensional even Banach spaces. Recall that the space $X_{GM}$ is the real version of the HI space of Gowers and Maurey \cite{GM1}, on which every operator is of the form $\lambda Id + S$, and therefore does not admit complex structure. The space $X(\C)$  is a HI space constructed in \cite{F} and such that
the algebra ${\cal L}(X(\C))$ may be decomposed as $\C \oplus {\cal S}(X(\C))$, and $X(\HH)$, also HI, is a quaternionic version of $X(\C)$. It is proved in \cite{F} that $X(\HH)$ admits a unique complex structure, while 
$X(\C)$ admits exactly two complex structures. Finally $X(\DD)$ is a space with an unconditional basis on which every operator is a strictly singular perturbation of a $2$-block diagonal operator and which also admits a unique complex structure. We refer to $\cite{F}$ for details.

\begin{theo}\label{existseven} The spaces $X(\C)$, $X(\HH)$, $X_{GM}^{2n}$ for $n \in \N$, and $X(\DD)$ are even. \end{theo}

\pf For $X(\C)$, $X(\HH)$ and $X_{GM}^{2n}$ this is due to the fact that 
the ideal of strictly singular operators has the lifting property, because
there is a natural subalgebra supplementing the ideal of strictly singular
operators in each case, see \cite{F}, so Lemma \ref{subalgebra} applies. For $X(\DD)$ note that it is proved in \cite{F} that
$X(\DD)$ admits a complex structure and that there is a unique conjugation class
in $\tilde{{\cal I}}_{\cal S}(X)$, therefore the induced injection $\tilde{\pi}_{\cal S}$ is
necessarily surjective, i.e. v) is verified. \pff

\

Before giving some more examples, let us note two open problems about even spaces. The first problem is quite simple to formulate.

\begin{ques} Is the direct sum of two even Banach spaces necessarily even?
\end{ques}

We obtain a positive answer when the spaces are assumed to be essentially incomparable. Recall that two infinite dimensional spaces $Y$ and $Z$ are essentially incomparable when every operator from $Y$ to $Z$ is inessential \cite{Gz}. More details about this notion may be found in the last section of this article. 

\begin{prop}\label{evenpluseven} The direct sum of two infinite dimensional even Banach spaces which are essentially incomparable is even. \end{prop}

\pf Let $X$ and $Y$ be infinite dimensional even, and essentially incomparable. Clearly $X \oplus Y$ admits a complex structure. Assume $X \oplus (Y \oplus \R)$ admits a complex structure and look for a contradiction.
Let $T=\begin{pmatrix} T_1 & S' \\ S & T_2 \end{pmatrix} \in {\cal L}(X \oplus (Y \oplus \R))$ be such that $T^2=-Id$. Since $X$ and $Y \oplus \R$ are essentially incomparable, $S$ and $S'$ are inessential. Furthermore
$T_1^2+S'S=-Id_X$, i.e. $\tilde{T_1} \in \tilde{{\cal I}}(X)$. Since $X$ is even,  there exists an inessential operator $s_1$ on $X$ such that
$T_1+s_1=J_1$ with $J_1^2=-Id_X$ (Theorem \ref{resume}).
Likewise $T_2^2+SS'=-Id_{Y \oplus \R}$ therefore since
$Y \oplus \R$ is odd, there exists an inessential operator $s_2$ on $Y \oplus \R$ such that $T_2+s_2=\begin{pmatrix} J_2 & 0 \\ 0 & 1 \end{pmatrix}$, $J_2 \in {\cal L}(Y)$ with
$J_2^2=-Id_Y$. 
Therefore there exists $S_0=\begin{pmatrix} s_1 & -S' \\ -S & s_2 \end{pmatrix}$ inessential on $X \oplus Y \oplus \R$ such that
$$T+S_0=\begin{pmatrix} J_1 & 0 & 0 \\ 0 & J_2 & 0 \\ 0 & 0 & 1 \end{pmatrix}.$$ 
Since $T^2=-Id$ this contradicts Proposition \ref{exclusion}.\pff

\

For the second problem, we note that Banach spaces which are not isomorphic to their hyperplanes  may be classified in four categories: even spaces, odd spaces, spaces such that neither the whole space neither hyperplanes admit a complex structure, and spaces such that both the whole space and any hyperplane admit complex structure.
While we have just produced examples of the first and the second category, and the space $X_{GM}$ belongs to the third, no examples are known which belong to the fourth.

\begin{ques}\label{q} Does there exist a real Banach space $X$ which is not isomorphic to its hyperplanes and such that both $X$ and its hyperplanes admit complex structure?
\end{ques}

\

We now use some spaces constructed by Plebanek \cite{P} to give $C(K)$ examples of even and odd Banach spaces. Similar $C(K)$ spaces were first constructed by P. Koszmider \cite{Kosz} under the Continuum Hypothesis.
Let $K$ be one of the two infinite, separable, compact Hausdorff spaces defined in \cite{P}. Every operator on $C(K)$ is of the form $g.Id+S$ where $g \in C(K)$ (therefore $g.Id$ denotes the multiplication by $g$) and $S$ is strictly singular (or equivalently weakly compact).
The first space is connected, and we shall indicate where our proofs simplify due to this additional property.

The space $K \cup K$ will denote the space which is the topological union of two copies of $K$ (i.e. open sets are unions of open sets of each copy), while $K \cup_0 K$ denotes the amalgamation of two copies of $K$ in some point $0$ (open sets are unions of open sets of each copy either both containing $0$ or neither containing $0$). Both are separable compact Hausdorff spaces.

\begin{theo}\label{CK} The space $C(K \cup K)$ is even and admits a unique complex structure, and the space $C(K \cup_0 K)$ is odd. \end{theo}

\pf The space $C(K \cup K)$ identifies isomorphically with $C(K) \oplus C(K)$ and
$C(K \cup_0 K)$ identifies with the quotient $C(K)^2/Y$, where $Y=\{(f,g): f(x_0)=g(x_0)\}$ for some fixed $x_0 \in K$; therefore $C(K \cup_0 K)$ is isomorphic to a hyperplane of
$C(K \cup K)$ and it is enough to prove that $C(K) \oplus C(K)$ is even with unique complex structure.

Write $X=C(K)$ and let $Is(K)$ be the set of isolated points of $K$.
 Let ${\cal N} \subset C(K)$ be the closed ideal of {\em almost null} functions, i.e. $g \in {\cal N}$ iff
$g$
vanishes on $K\setminus Is(K)$  and converges to $0$ on $Is(K)$ (i.e. for any $\epsilon>0$, $|g(x)| \leq \epsilon$ for all $x \in Is(K)$ except a finite number of points).
We observe the following fact:
if an operator on $X$ of the form $g. Id$ is strictly singular, then $g$ belongs to ${\cal N}$.
Indeed if $g(x) \neq 0$ for some non-isolated point $x$ then $|g| \geq \alpha >0$ on some
infinite subset $L$ of $K$ containing $x$, and the restriction of $g Id$ to 
the space of functions with support included in $L$ is an isomorphism.
Likewise if $|g| \geq \alpha>0$ on some infinite subset of $Isol(K)$ then the corresponding restriction of $g Id$ is an isomorphism.

In the case where $K$ is connected, we have simply that $g. Id$ is never strictly singular unless $g=0$.

Since $X \oplus X$ has a canonical complex structure, to prove that it is even with unique complex structure, it is enough by Theorem \ref{resume} to prove that the group $G_0:=({\cal L}(X \oplus X)/{\cal S}(X \oplus X))_0$ has a unique conjugation class of elements of square $-Id$.

\

We observe that the group $GL_2(C(K))$ admits a unique class of conjugation of elements of square $-I$, where $I=\begin{pmatrix} 1 & 0 \\ 0 & 1 \end{pmatrix}$, i.e.
 that whenever $M=\begin{pmatrix} f_1 & f_2 \\ f_3 & f_4 \end{pmatrix}$ satisfies $M^2=-I$, then it is conjugate to the canonical element $J=\begin{pmatrix} 0 & 1 \\ -1 & 0 \end{pmatrix}$.
Indeed from $M^2=-I$
 we deduce easily 
$f_1=-f_4$ and $f_1^2+f_2f_3=-1$.
Note that $f_2$ never takes the value $0$.
Let then $P=\begin{pmatrix} 1 & 0 \\ f_1 & f_2 \end{pmatrix}$
and let
$Q=\begin{pmatrix} 1 & 0 \\ -f_1/f_2 & 1/f_2 \end{pmatrix}.$
 It is routine to check that $Q=P^{-1}$ and that
$QJP=M$.

 When $K$ is connected then ${\cal N}=\{0\}$, therefore
${\cal L}(X) \simeq C(K) \oplus {\cal S}(X)$
and ${\cal L}(X \oplus X)=M_2(C(K)) \oplus {\cal S}(X \oplus X).$
So ${\cal S}(X \oplus X)$ has the lifting property by Lemma \ref{subalgebra}, and $G_0$ easily identifies with $GL_2(C(K))$, and we therefore deduce that there is a unique $G_0$-conjugacy class of elements of square $-\tilde{Id}$.

\

 The general case is more complicated. 
Assume $\tilde{T} \in G_0$ satisfies 
$\tilde{T}^2=-\tilde{Id}$. Then $T^2+Id$
is strictly singular, and up to a strictly singular perturbation we may
assume that $T=M.Id$ with $M \in M_2(C(K))$. 
Therefore $(M^2+I).Id$ is strictly singular and $M^2+I \in {\cal N}$.
It follows that
$M^2=-I+n$ with $\||n(x)\|| \leq 1/2$ except for $x \in F$, $F$ a finite subset of $Is(K)$ (here $\||.\||$ denotes the operator norm on $M_2(\R) \simeq {\cal L}(\R \oplus_2 \R)$). Note that $M(x)$ and $n(x)$ commute for all $x \in K$.

For $x \in K\setminus F$, let $n'(x)=M(x)\sum_{k \geq 1}b_k n(x)^k$ where
$(1-z)^{-1/2}=1+\sum_{k \geq 1} b_k z^k$ for $|z| <1$, and let
$M'(x)=(M+n')(x)$. Then it is easy to check that
$(M'(x))^2=-Id_{\R^2}$ and therefore there exist $P(x), Q(x)$ in $M_2(\R)$ given by formulas from the coefficients of $M'(x)$ which are explicited above in the connected case, such that $P(x)Q(x)=Id_{\R^2}$ and $Q(x)jP(x)=M'(x)$, where 
$j:=\begin{pmatrix} 0 & 1 \\ -1 & 0 \end{pmatrix}$.

For $x \in F$ we let $M'(x)=(M+n')(x)=j$ and $P(x)=Q(x)=Id_{\R^2}$. 
Note that since the points of $F$ are isolated and by uniform convergence of $n'$ and the formulas giving $P$ and $Q$, the matrices
$P$, $Q$ and $n'$ define elements of $M_2(C(K))$.
Actually, since $M'(x)$ is invertible for all $x$ and by the explicit formulas for the inverses, we deduce that $M' \in GL_2(C(K))$ and
$M'.Id \in {\cal GL}(X \oplus X)$. Likewise $P$ and $Q$ belong to
$GL_2(C(K))$ and $P. Id$, $Q.Id$ belong to
${\cal GL}(X \oplus X)$.

It is now enough to prove that $S=n'.Id$ is a strictly singular operator on $X \oplus X$.  Then
the relation $QJP=M'$ will imply a $G_0$-conjugacy relation between 
$\tilde{J.Id}$ and $\tilde{M'.Id}=\tilde{M.Id}+\tilde{S}=\tilde{T}$, as desired.

For $L \subset K$, $C_L(K)$ denotes the space of functions of support included in $L$. Let $v$ denote the rank $|F|$ projection onto $C_F(K)$
associated to the decomposition $C(K)=C_F(K) \oplus C_{K\setminus F}(K)$,
and $w=Id-v$. Let $V$ on $X \oplus X$ be defined by $V= \begin{pmatrix} v & 0 \\ 0 & v \end{pmatrix}$
and let $W=\begin{pmatrix} w & 0 \\ 0 & w \end{pmatrix}$.

It is easy to check that $S=\sum_{k \geq 1}b_k(n^k.Id)W+V(n'.Id)$, (just  compute $S(f_1,f_2)(x)$ for all $(f_1,f_2) \in X \oplus X$ and $x \in K$). Since $n^k.Id=(n.Id)^k$ is strictly singular for all $k$ and $V$ has finite rank, it follows that $S$ is strictly singular. \pff

\

To conclude this section we note an open question. Recall that there exist spaces with exactly $n$ complex structures, for any $n \in \N^*$ \cite{F}, and spaces with exactly $2^{\omega}$ complex structures \cite{An}.

\begin{ques} Does there exist a Banach space with exactly $\omega$ complex structures?
\end{ques}

\section{Essentially incomparable complex structures}

We recall that two infinite dimensional spaces $Y$ and $Z$ are said to be {\em essentially incomparable} if
any bounded operator $S \in {\cal L}(Y,Z)$ is inessential, i.e., if
$Id_Y-VS$ is a Fredholm operator (necessarily of index $0$) for all operators $V \in {\cal L}(Z,Y)$. Essentially
incomparable spaces were studied by M. Gonz\'alez in \cite{Gz}; it is clear that the notion of
essential incomparability generalizes the notion of total incomparability.
We also recall that $Y$ and $Z$ are {\em projection totally incomparable} if no infinite dimensional complemented subspace of $Y$ is isomorphic to a complemented subspace of $Y$. Essentially incomparable spaces are in particular projection totally incomparable.

In this section we show how some results of \cite{F} about totally
incomparable complex structures extend to essentially incomparable structures. Interestingly, our more general proof turns out to be much simpler than the original one.

As was noted in \cite{F}, whenever $T^{2}=U^{2}=-Id$, it follows that
$$(T+U)T=U(T+U),$$ which means that $T+U$ is $\C$-linear from
$X^T$ into $X^U$. The similar result
holds for $T-U$ between $X^T$ and $X^{-U}$.
\

\begin{lemm}\label{strange}
Let $X$ be a real Banach space, $T,U \in {\cal I}(X)$.  If $T+U$ is inessential from $X^T$ into $X^{U}$, then $X^T$ is
isomorphic to
 $X^{-U}$. 
\end{lemm}

\pf 
 Since $T+U$ is inessential as an operator from $X^T$ into $X^U$, and as $T+U$ is also linear from $X^U$ into $X^T$, it follows by definition of inessential operators that $Id+\lambda(T+U)^2 \in {\cal L}(X^T)$ is Fredholm with index $0$ whenever $\lambda$ is real. Taking $\lambda=1/4$, we obtain that
$$4Id+(T+U)^2=2Id+TU+UT=-(T-U)^2.$$
Therefore $(T-U)^2$ is Fredholm with index $0$ as an operator on $X^T$, and therefore as an operator on $X$. It follows that
$T-U$ is Fredholm with index $0$ as an operator of ${\cal L}(X)$ and therefore as an operator of ${\cal L}(X^T,X^{-U})$, hence
$X^T$ and $X^{-U}$ are $\C$-linearly isomorphic.
\pff

\

It was proved in \cite{F} that two totally incomparable complex structures on a real space must be conjugate and both saturated with HI subspaces. We show:

\begin{theo}\label{ttt}
Let $X$ be a real Banach space with two essentially incomparable complex structures. Then
these complex structures are conjugate up to isomorphism and do not contain a
complemented subspace with an unconditional basis.
\end{theo}

\pf 
Assume $X^T$
is essentially  incomparable with
$X^U$. Then $T+U$ is inessential from $X^T$ into $X^U$ and
by Proposition
\ref{strange}, $X^U$ is isomorphic to $X^{-T}$.

If $Y$ is a $\C$-linear  complemented subspace of $X^{T}$ with an
unconditional basis, then $\overline{Y}$ is complemented in $X^{-T}$ and
the coordinatewise conjugation map $\alpha$ associated to the unconditional basis is an isomorphism from $Y$
onto $\overline{Y}$. Therefore $X^T$ and $X^{-T}$ are not
projection totally incomparable,
contradicting the essential incomparability of
$X^T$ with $X^U \simeq X^{-T}$. 
 \pff

\

Note that Proposition \ref{exam} will prove that one cannot hope to improve Theorem \ref{ttt} to obtain HI-saturated in its conclusion, as in the case of totally incomparable complex structures.

\begin{coro}\label{ccc}
There cannot exist more than two mutually essentially incomparable complex structures on
a Banach space.
\end{coro}

Recall that two Banach spaces are said to be {\em nearly isomorphic} (or sometimes essentially isomorphic) if one is isomorphic to a finite-codimensional subspace of the other. Equivalently this means that there exists a Fredholm operator acting between them.

In the next proposition, we consider properties of complex structures which are generalization of the $\lambda Id+S$-property. Note that each of these properties implies that there do not exist non-trivial complemented subspaces. We first state a lemma whose proof was given to us by M. Gonz\'alez.

\begin{lemm}\label{gon} Let $X$ be an infinite dimensional real or complex Banach space such that every operator is either Fredholm or inessential. Then every Fredholm operator on $X$ has index $0$ and ${\cal L}(X)/{{\cal I} n}(X)$ is a division algebra. \end{lemm}

\pf If in the above conditions, $T$ were Fredholm with nonzero index, then by the continuity of the index, $K:=sT+(1-s)Id$ would be inessential for some $s$ with $0<s<1$. Thus $T:=\frac{s-1}{s} Id+\frac{1}{s} K$; hence $T$ has index $0$, a contradiction. Moreover, the only noninvertible element in  ${\cal L}(X)/{{\cal I}n}(X)$ is $0$; hence
${\cal L}(X)/{{\cal I}n}(X)$ is a division algebra. \pff

\begin{prop}\label{1to5}
Let $X$ be a real Banach space, and $T \in {\cal I}(X)$.
\begin{itemize}

\item i) If every operator on $X^T$ is either inessential or Fredholm, then either $X^T$ is the only complex structure on $X$,
or $X^T$ and $X^{-T}$ are the only two complex structures on $X$ and they are not nearly isomorphic.

\item ii) If every operator on $X^T$ is either strictly singular or Fredholm, then either $X^T$ is the only complex structure on $X$,
or $X^T$ and $X^{-T}$ are the only two complex structures on $X$ and  neither one  embeds  into the other.

\end{itemize}
\end{prop}

\pf i) Let $U$ generate a complex structure on $X$. We use the relation
$$(T+U)^2+(T-U)^2=-4 Id.$$
If $(T-U)^2$, which is an operator on $X^T$, is inessential, then
$(T+U)^2$ is Fredholm with index $0$ as a perturbation of $-4 Id$, therefore $T+U$ is Fredholm with index $0$ and there exists an isomorphism from $X^T$ onto $X^U$.

Therefore if there exists some $U$ generating a  complex structure non
isomorphic to $X^T$, then $(T-U)^2$ is not inessential, and by the property of operators on $X^T$, $(T-U)^2$ is Fredholm, with index $0$ by Lemma \ref{gon}, and therefore $T-U$ as well, so $X^U$ is isomorphic to $X^{-T}$.
We deduce that $X^T$ and $X^{-T}$ are the only complex structures on $X$. To see that they are not nearly isomorphic, note that if a map $\alpha$ is Fredholm from $X^T$ into $X^{-T}$, then $\alpha^2$ is Fredholm on $X^T$, hence it is Fredholm with index $0$, and $\alpha$ is Fredholm with index $0$. Therefore $X^T$ is isomorphic to $X^{-T}$, contradicting our initial assumption.

ii) Applying i) we see that any complex structure on $X$ is either isomorphic to $X^T$ or $X^{-T}$, and if $\alpha$ embeds $X^T$ into $X^{-T}$, then $\alpha^2$ embeds $X^T$ into itself, hence it is not strictly singular and so it is Fredholm, with index $0$, and $\alpha$ is Fredholm with index $0$, which implies that $X^T$ is isomorphic to $X^{-T}$.
 \pff

\

In \cite{AM} S. Argyros and A. Manoussakis constructed a real space  which is unconditionnally saturated yet has the $\lambda Id+S$ property. Although their result is stated in the real case, no specific property of the reals is used in their definition and proofs, and so the complex version $X_{AM}$ of their space satisfies the complex version of the properties mentioned above.
We observe:

\begin{prop}\label{exam} Every $\R$-linear operator on the complex $X_{AM}$ is of the form $\lambda Id+S$, $\lambda \in \C$, $S$ strictly singular. It follows that the complex $X_{AM}$ seen as real admits exactly two complex structures, which are essentially incomparable but not totally incomparable.
\end{prop}

\pf Denote $X_{AM}$ the complex version of the space of Argyros and Manoussakis, and let $X$ be $X_{AM}$ seen as real.
Lemma 4.17 from \cite{AM} in its complex version states that every $\C$-linear operator $T$ on $X_{AM}$ satisfies $\lim d(Te_n,\C e_n)=0$, where $(e_n)_n$ is the canonical complex basis of $X_{AM}$. A look at their proof shows that actually only the $\R$-linearity of $T$ is required. Then if $T$ is $\R$-linear on $X$ we deduce easily that there exists $\lambda \in \C$ such that for any $\epsilon>0$, there exists $M$ an infinite subset of $\N$ such that $(T-\lambda Id)_{|[e_n, n \in M]}$ is of norm at most $\epsilon$.
Here $[e_n, n \in M]$ denotes the real subspace $\R$-linearly generated
by $(e_n)_{n \in M}$.
 
Proposition 4.16 in \cite{AM} states that for any infinite
subset $M$ of $\N$, any $(y_k)$ a normalized block-sequence of $(e_n)$,
the distance
$d(S_{[e_n,n \in M]},S_{[y_k,k \in \N]})$ between the respective unit spheres
of the complex (i.e. $\C$-linearly generated) block-subspaces $[e_n, n \in M]$ and $[y_k, k \in \N]$ is $0$.
 A look at the proof shows that one can obtain this using only the real block-subspaces which are $\R$-linearly generated by $(e_n)_{n \in M}$ and $(y_k)_{k \in \N}$ (in particular note that \cite{AM} Lemma 4.12 used in the proof only uses $\R$-linear combinations of the $(e_n)$'s and the $(y_k)$'s).

Combining the facts of the first and the second paragraph, we deduce that when $T$ is $\R$-linear on $X$, there exists $\lambda \in \C$ such that for all $\epsilon>0$, for any 
complex block-subspace $(y_k)_k$ of $(e_n)_n$, there exists a unit vector
$x$ in the $\R$-linear span of $(y_k)_k$ such that
$(T-\lambda Id)x$ is of norm less than $2\epsilon$. We deduce easily that
every $\R$-linear operator on $X$ is of the form $\lambda Id+S$, $\lambda \in \C$, $S$ strictly singular. 

The complex structure properties now follow easily.
By Proposition \ref{1to5} ii) we already know that either
$X$ has unique complex structure, or exactly two which are $X^J$ and $X^{-J}$ (where $J \in {\cal L}(X)$ is defined by $Jx=ix$). Note that whenever $\alpha$ is $\C$-linear
from $X^J$ into $X^{-J}$ then $\alpha(ix)=-i\alpha(x)$ for all $x \in X$. Since
$\alpha=\lambda.Id+S$, $\lambda \in \C$ and $S$ $\R$-strictly singular, we deduce that $\lambda=0$ and $\alpha=S$. The operator $S$ is in particular $\C$-strictly singular;
so ${\cal L}(X^J,X^{-J})={\cal S}(X^J,X^{-J})$ and $X^J$ and $X^{-J}$ are essentially incomparable. On the other hand 
$X^J$ is the complex version of $X_{AM}$ and so is unconditionally saturated, therefore it is not HI saturated, and so not totally incomparable with its conjugate $X^{-J}$, by \cite{F} Corollary 23.\pff

\

We end the paper with two open questions in the direction of further generalizing the above results. 
Recall that
two spaces are projection totally incomparable if no infinite dimensional complemented subspace of one is isomorphic to a complemented subspace of the other, and that essentially incomparable spaces are in
particular projection totally incomparable \cite{Gz}.

\begin{ques} If two complex structures on a real space $X$ are projection totally incomparable, must they be conjugate?
\end{ques}

\begin{ques} Assume a complex space is projection totally incomparable with its conjugate, is it necessarily essentially incomparable with it?
\end{ques}

\paragraph{Acknowledgements} The first author thanks B. Maurey for a useful discussion which led the authors to Proposition  \ref{exclusion}. Both authors thank M. Gonz\'alez for his comments on a first version of this paper, and P. Koszmider for information about the
$C(K)$ examples used in Theorem 
\ref{CK}.

\

\

Valentin Ferenczi,

Institut de Math\'ematiques de Jussieu,

Universit\'e Pierre et Marie Curie - Paris 6,

Projet Analyse Fonctionnelle, Bo\^\i te 186,

4, place Jussieu, 75252 Paris Cedex 05,

France.

\

E-mail: ferenczi@ccr.jussieu.fr, ferenczi@ime.usp.br

\

\

El\'oi Medina Galego,

Departamento de Matem\'atica,

Instituto de Matem\'atica e Estat\' \i stica,

Universidade de S\~ao Paulo.

05311-970 S\~ao Paulo, SP,

Brasil.

\

E-mail: eloi@ime.usp.br.


\begin{thebibliography}{Wwww}

\bibitem{An} R. Anisca, {\em Subspaces of $L_p$ with more than one complex structure},
Proc. Amer. Math. Soc.  {\bf 131}  (2003),  no. 9, 2819--2829.


 \bibitem{AM} S. Argyros and A. Manoussakis, 
{\em An indecomposable and unconditionally saturated Banach space},
Studia Math. {\bf 159} (2003), no. 1, 1--32.


\bibitem{B} J. Bourgain, {\em Real isomorphic complex Banach
spaces need not be complex isomorphic},  Proc. Amer. Math. Soc. 
{\bf 96}  (1986),  no. 2, 221--226.


 \bibitem{DS}  N. Dunford and J. T. Schwarz. {\em Linear Operators, Part I: General theory}, New York 1958.


\bibitem{F} V. Ferenczi, {\em Uniqueness of complex structure and real hereditarily indecomposable Banach spaces}, Advances in Math., to appear. 
 
\bibitem{Gon}  M. Gonz\'alez, {\em Banach spaces with small Calkin algebras}, Banach Center Publ. to appear.

\bibitem{Gon2} M. Gonz\'alez and J. M. Herrera. {\em Decompositions for real Banach spaces with small spaces of operators}, preprint.

\bibitem{Gz} M. Gonz\'alez, {\em On essentially incomparable Banach spaces},  Math. Z.  {\bf 215}  (1994),  no. 4, 621--629.
 
\bibitem{GM1} W.T. Gowers and B. Maurey, {\em The unconditional basic
sequence problem}, J. Amer.
Math. Soc. {\bf 6} (1993), 4, 851--874.

\bibitem{K} N. J. Kalton, {\em An elementary example of a Banach space not isomorphic to its
complex conjugate},  Canad. Math. Bull.  {\bf 38} (1995),  no. 2, 218--222.

\bibitem{Kosz} P. Koszmider, {\em Banach spaces of continuous functions with few operators}, Math. Ann. {\bf 330} (2004), 151--183.

\bibitem{LT} J. Lindenstrauss and L. Tzafriri, {\em Classical Banach
spaces}, Springer-Verlag, New York, Heidelberg, Berlin (1979).
  

\bibitem{P} G. Plebanek, {\em A construction of a Banach space $C(K)$ with few operators}, Topology and its applications {\bf 143} (2004), 217--239.



\end{thebibliography}
\end{document}